\documentclass[12pt]{amsart}

\usepackage{yhmath}

\usepackage{amssymb,amsmath,amsfonts,amscd,amsthm,enumerate}
\usepackage{graphicx,color}
\graphicspath{{picturesfolder/}}
\usepackage{wrapfig}
\usepackage{multirow}
\usepackage{textcomp}

\usepackage{todonotes}

\usepackage{epsfig,subfigure}
\usepackage{todonotes}

\usepackage{lineno,xcolor}

\usepackage[all]{xy}
\usepackage{color}

\newcommand{\R}{{\Bbb R}}
\newcommand{\C}{{\Bbb C}}

\newcommand{\so}{\Rightarrow}
\newcommand{\e}{\varepsilon}

\newcommand{\inv}{^{-1}}

\newcommand{\del}{\partial}

\newtheorem{theorem}{Theorem}

\newtheorem{lem}[equation]{Lemma}

\title{Hyperbolic Shirts fit a 4-body problem}

\author[Connor Jackman]{Connor Jackman}
\address[Jackman]{Mathematics Department, University of California,
4111 McHenry
Santa Cruz, CA 95064, USA}
\email{cfjackma@ucsc.edu}

\author[Josu\'e Mel\'endez]{Josu\'e Mel\'endez}
\address[Mel\'endez]{Departamento de Matem\'aticas, UAM-Iztapalapa, 09340 M\'exico City, M\'exico}
\email{jms@xanum.uam.mx}

\begin{document}

\date{\today}
  
\begin{abstract}  Consider the equal mass planar $4$-body problem with a potential corresponding to an inverse \textit{cube} force. The Jacobi-Maupertuis principle reparametrizes the dynamics as geodesics of a certain metric. We examine the curvature of this geodesic flow in the reduced space on the collinear and parallelogram invariant surfaces and derive some dynamical consequences. This proves a numerical conjecture of \cite{Nopants}.
\end{abstract}

\maketitle



\section{Introduction}
\label{intro}

In \cite{Nopants} it was shown that the reduced Jacobi-Maupertuis metric for the 4-body and 3-body strong force problems have different curvature properties. The 3-body metric was shown in \cite{Pants} to be negatively curved (except at two points) leading to interesting dynamical consequences. Yet for $N\ge 4$ we found in \cite{Nopants} that this negativity does not hold, the sectional curvatures of the reduced Jacobi-Maupertuis metric have mixed signs.

Where then is the curvature negative? Are there any dynamical consequences? Here we will show that for $N=4$ the reduced Jacobi-Maupertuis metric is negatively curved (except at two points) over the invariant surfaces corresponding to parallelogram or collinear configurations and derive some dynamical consequences.

\section{Notations and Outline of Results}
\label{sec2}

We consider the equal mass ($m_i=1$) `strong force' planar 4-body problem with configuration space $$q=(q_1,q_2,q_3,q_4)\in \C^4\backslash\Delta$$ where $\Delta=\{ q: q_i=q_j\text{ for some } i\ne j\}$ are the collisions. The dynamics are then governed by the Hamiltonian flow of $$H=\frac12|p|^2-U,\quad \omega=\Re (d\overline p\wedge dq)$$ where $|\cdot |$ is the euclidean norm and $$U=\sum_{i<j} |q_i-q_j|^{-2}$$ is a `strong force' potential.

As is standard, due to the translation invariance of $H$, we consider solutions in center of mass zero coordinates that is in the invariant $$CM_0:=\{ q: \sum q_i=0\}\backslash \Delta.$$

Due to the homogeneity of degree $-2$ of the strong force the virial identity reads as $$\ddot I=4H$$ where $I=|q|^2$. Hence the energy zero case is the interesting case with regard to finding periodic motions.

The Jacobi-Maupertuis principle (see \cite{ACM} pg. 247) states that the solutions at this fixed energy level $H=0$ are upto reparametrazation geodesics of the metric \begin{equation}
    Uds_{Eucl}^2
    \label{JM}
\end{equation} restricted to $CM_0\backslash\Delta$.

Again, due to $U$'s homogeneity of degree $-2$, the metric (\ref{JM}) is invariant under complex multiplication so the Hopf map: $$\pi_{Hopf}:CM_0\backslash\Delta\to \C P^2\backslash P\Delta$$ pushes down (\ref{JM}) by submersion to obtain a metric $ds_{JM}^2$ on the quotient space which we call the reduced Jacobi-Maupertuis metric (or JM-metric) on the shape space. See \cite{Shape}, \cite{Pants}, \cite{Nopants} for more on the shape space.\\

{\sc Remark:} In \cite{Pants} the reduced strong force Jacobi-Maupertuis metrics are shown to be complete for all $N$.

{\sc Remark:} Pushing down by submersion here is equivalent to imposing the conditions $\dot I=0$ and angular momentum, $J$, zero ($J=\sum_j ds_{Eucl}^2(iq_j,\dot q_j)$) on the solutions.\\

Here we show that the curvature of this $ds_{JM}^2$ metric over the totally geodesic surfaces $$\mathcal{P}=\pi_{Hopf}(\text{Fix}\{ q\mapsto -(q_3,q_4,q_1,q_2)\})\cong S^2\backslash \{ 4\,\text{points}\}$$ and $$\mathcal{C}=\pi_{Hopf}(\text{Fix}\{q\mapsto \overline q\})\cong \R P^2\backslash RP\Delta$$ are non-positive:

\begin{theorem}  The curvature of $ds_{JM}^2$ over $\mathcal{P}$ is negative except at two points corresponding to the square configurations, where it is zero.
\label{P}
\end{theorem}

\begin{theorem}  The curvature of $ds_{JM}^2$ over $\mathcal{C}$ is strictly negative.
\label{C}
\end{theorem}

{\sc Remark} The set $\mathcal{P}$ consists of configurations that are parallelograms (there are two other such parallelogram subspaces corresponding to different choices of diagonals), and the set $\mathcal{C}$ consists of configurations that are collinear.

\begin{figure}[ht!]
\centering
\includegraphics[height=75mm,width=90mm] {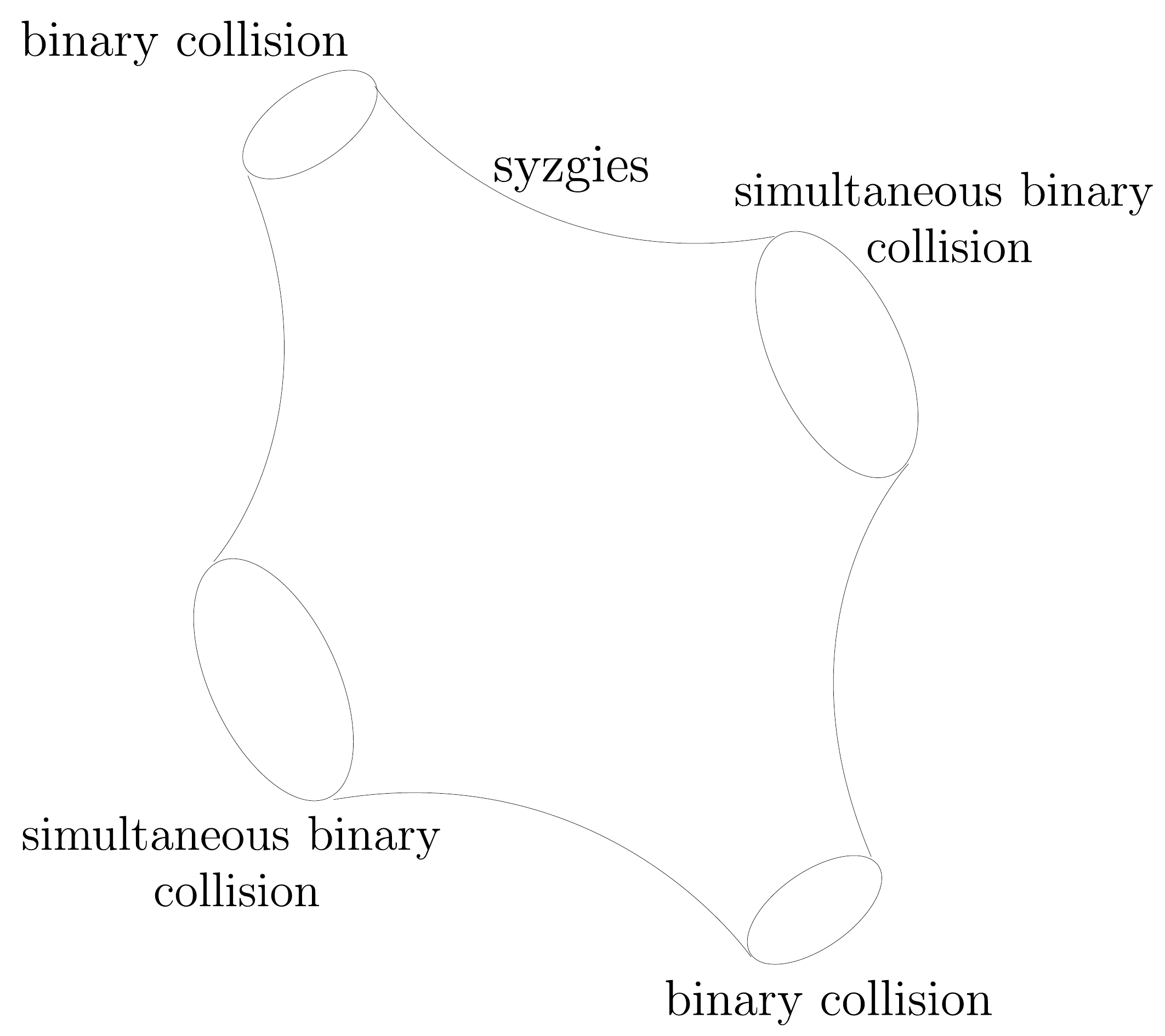}
\caption{~~The shirt.}
\label{camisa}
\end{figure}

{\sc Remark} We follow the techniques used in \cite{AGMR} to prove Theorem \ref{P}. Note that the surface $\mathcal{P}$ is topologically equivalent to a $2$-sphere minus four points, or a shirt. The punctures of the $2$-sphere are due to two binary collisions ($q_1=q_3$ and $q_{2}=q_4$) and two simultaneous binary collisions ($q_1=q_4, q_2=q_3$, and  $q_1=q_2, q_3=q_4$). See Figure \ref{camisa} for a depiction of $\mathcal{P}$.

\begin{figure}
\centering
\begin{tikzpicture}[scale=0.70]
\draw (-6,0)-- (-4.5,0) node[below] {$m_1$} -- (-1.5,0) node[below] {$m_2$} -- (1.5,0) node[below] {$m_4$} -- (4.5,0) node[below] {$m_3$}--(6,0);
\shade[shading=ball, ball color=black] (-4.5,0) circle (.16);
\shade[shading=ball, ball color=black] (-1.5,0) circle (.16);
\shade[shading=ball, ball color=black] (1.5,0) circle (.16);
\shade[shading=ball, ball color=black] (4.5,0) circle (.16);
\draw (-6,-2)-- (-4.5,-2) node[below] {$m_1$} -- (-1.5,-2) node[below] {$m_4$} -- (1.5,-2) node[below] {$m_2$} -- (4.5,-2) node[below] {$m_3$}--(6,-2);
\shade[shading=ball, ball color=black] (-4.5,-2) circle (.16);
\shade[shading=ball, ball color=black] (-1.5,-2) circle (.16);
\shade[shading=ball, ball color=black] (1.5,-2) circle (.16);
\shade[shading=ball, ball color=black] (4.5,-2) circle (.16);
\draw (-6,-4)-- (-4.5,-4) node[below] {$m_4$} -- (-1.5,-4) node[below] {$m_1$} -- (1.5,-4) node[below] {$m_3$} -- (4.5,-4) node[below] {$m_2$}--(6,-4);
\shade[shading=ball, ball color=black] (-4.5,-4) circle (.16);
\shade[shading=ball, ball color=black] (-1.5,-4) circle (.16);
\shade[shading=ball, ball color=black] (1.5,-4) circle (.16);
\shade[shading=ball, ball color=black] (4.5,-4) circle (.16);
\draw (-6,-6)-- (-4.5,-6) node[below] {$m_4$} -- (-1.5,-6) node[below] {$m_3$} -- (1.5,-6) node[below] {$m_1$} -- (4.5,-6) node[below] {$m_2$}--(6,-6);
\shade[shading=ball, ball color=black] (-4.5,-6) circle (.16);
\shade[shading=ball, ball color=black] (-1.5,-6) circle (.16);
\shade[shading=ball, ball color=black] (1.5,-6) circle (.16);
\shade[shading=ball, ball color=black] (4.5,-6) circle (.16);
\end{tikzpicture}
\caption{~~The four syzygy types on $\mathcal{P}$.}
\label{syzygy}
\end{figure}
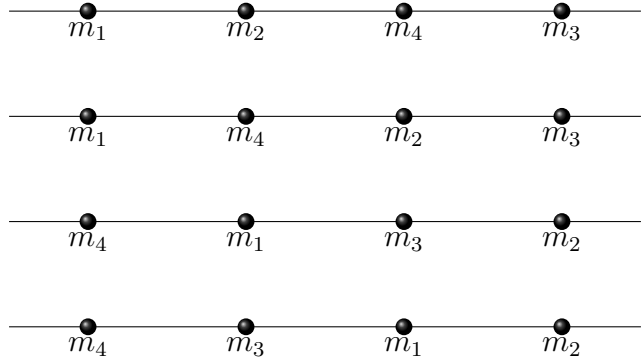

A syzygy is a collinear configuration on the equator of $\mathcal{P}$. See Figure \ref{syzygy}. The syzygy map takes a solution and lists its syzygy types in temporal order. Deleting the stutters makes this map well defined on free homotopy clases (see \cite{Pants} \S 2 pgs. 3-6).

Theorem \ref{P} implies the following: (in the same way as \cite{Pants})

\begin{theorem}
The syzygy map from bounded parallelogram solutions to syzygy sequences is a bijection between the set of collision-free solutions, modulo symmetry and time-translation, and the
set of bi-infinite non-stuttering collision-free syzygy sequences, modulo shift.
\label{J}
\end{theorem}

{\sc Remark} The surface $\mathcal{C}$ consists of 12 invariant open disks, each disk corresponding to an ordering of the masses on the line (mod reflections).\\

If one restricts attention to the closure of one of these invariant disks, $T\subset\mathcal{C}$ (such as $\{q: q_1\le q_2\le q_3\le q_4\}$), we prove here that:

\begin{theorem}
Let $p, q\in\del T\backslash \{$ the two triple collision points $\}$ be two distinct points. Then there exists a unique solution, $\gamma(t)$ such that $\omega(\gamma)=p$ and $\alpha(\gamma)=q$.
\label{M}
\end{theorem}

{\sc Open Questions} Do analogous results above hold for unequal masses? For $N>4$? Does Theorem \ref{M} hold when $p$ and $q$ are triple collisions? Does Theorem \ref{J} hold for collision syzygy sequences? See Figures \ref{sol1}, \ref{sol2} for some collision orbits with matching syzygy sequences (giving some numerical evidence that the syzygy map is not 1-1 on collision orbits). Are there any solutions besides the rectangulars and rhomboids which have finite syzygy sequences? 

\section{Proof of Theorems}
\label{sec3}

\textbf{Proof of Thm. \ref{P}:}
Consider the symmetry of the problem: $q_1=-q_3$ and $q_2=-q_4$. In this case the (negative) potential function has
the following form
\begin{equation}\label{eq:U}
U=\frac{2}{r_{12}^2}+\frac{2}{r_{14}^2}+\frac{1}{r_{13}^2}+\frac{1}{r_{24}^2}.
\end{equation}
Next we take the quotient by rotations. As in \cite{Chen}, consider Jacobi's coordinates and the Hopf map:
\begin{eqnarray*}
(q_1,q_2,q_3,q_4) &\mapsto& (q_2-q_1,-q_2-q_1)=:(z_1,z_2) \in \mathbb{C}^2, \\
(z_1,z_2) &\mapsto& \left(\frac{1}{2}\left(\left|z_1\right|^2-\left|z_2\right|^2 \right),\bar{z}_1z_2\right)=:	(u_1,u_2,u_3) \in \mathbb{R}^3. \nonumber
\end{eqnarray*}
Using the above formulas, it is immediate to see that (\cite[Theorem 1]{Shape})
\[
I=|z_1|^2+|z_2|^2=2\sqrt{u_1^2+u_2^2+u_3^2}, \qquad I=\text{constant}.
\]

The following Lemma describes the relation between points of $\mathcal{P}$
and the mutual distances $r_{ij}$ (see \cite{Chen}):

\begin{lem}
\label{lem:distances}
 To $u = (u_1, u_2, u_3)\in\mathcal{P}$, we have
\begin{enumerate}
\item $r_{12}^2=|z_1|^2=\frac{I}{2}+u_1$,
\item $r_{14}^2=|z_2|^2=\frac{I}{2}-u_1$,
\item $r_{13}^2=|z_1+z_2|^2=I+2u_2$,
\item $r_{24}^2=|z_1-z_2|^2=I-2u_2$.
\end{enumerate}
\end{lem}

In particular, observe that $u_1 = 0$ if and only if the configuration is a rhombus, and
$u_2 = 0$ if and only if the configuration is a rectangle. So, we have that on the north or south pole the configuration is a square.

We will use the following Lemma in our computations:

\begin{lem}
\label{lem:crmetric}
Let a surface be endowed with conformally related metrics $ds^2$ and $d\overline{s}^2=Uds^2$. Then their curvatures
$K$ and $\overline{K}$ are related by
\[
U \overline{K}=K-\frac{1}{2} \Delta \log U
\]
where the Laplacian $\Delta$ is with respect to the $ds^2$ metric.
\end{lem}

From equation (\ref{eq:U}) and Lemma \ref{lem:distances}, the potential can be written as
\begin{equation}\label{eq:u}
U(u_1,u_2)=\frac{4}{I+2u_1}+\frac{4}{I-2u_1}+\frac{1}{I+2u_2}+\frac{1}{I-2u_2}, \qquad I=\text{constant}.
\end{equation}
Taking $I=2$, we have $u_1^2+u_2^2+u_3^2=1$. As in \cite{AGMR}, we use the stereographic projection from north pole $(0,0,1)$ given by
\begin{equation}\label{eq:proj}
x=\frac{u_1}{1-u_3}, \qquad y=\frac{u_2}{1-u_3},
\end{equation}
we obtain $U(x,y)=(x^2+y^2+1)f(x,y)$, where
\[
f(x,y)=\frac{2}{(x+1)^2+y^2}+\frac{2}{(x-1)^2+y^2}+\frac{1}{2(x^2+(y+1)^2)}+\frac{1}{2(x^2+(y-1)^2)}.
\]
We can also write the Jacobi-Maupertuis metric on the shape sphere as

\begin{equation}
ds^2_{JM}=\frac{4U(x,y)}{(x^2+y^2+1)^2}(dx^2+dy^2)=\frac{4f(x,y)}{x^2+y^2+1}(dx^2+dy^2).
\end{equation}
Write $\lambda(x,y)=\frac{4f(x,y)}{x^2+y^2+1}$. By Lemma \ref{lem:crmetric}, the Gaussian curvature of $ds^2_{JM}$ is given by
\begin{equation}\label{eq:K}
K(x,y)=-\frac{1}{2 \lambda(x,y)} \Delta \log \lambda(x,y).
\end{equation}
By a direct calculation we have
\begin{equation}  \label{eq:Dlambda}
\Delta  \log \lambda(x,y)=\frac{256 (x^2 + y^2)}{[5 + 5 x^4 - 6 y^2 + 5 y^4 + 2 x^2 (3 + 5 y^2)]^2}.
\end{equation}

Note that $\Delta \log \lambda(x,y)\geq 0$, and $\Delta \log \lambda(x,y)=0$ if and only if $x=y=0$. Therefore we conclude that  $K$
is negative except at the points corresponding to the square configurations.
\qed

\begin{figure}[htb]
\begin{center}
\includegraphics[width=8cm, height =5cm]{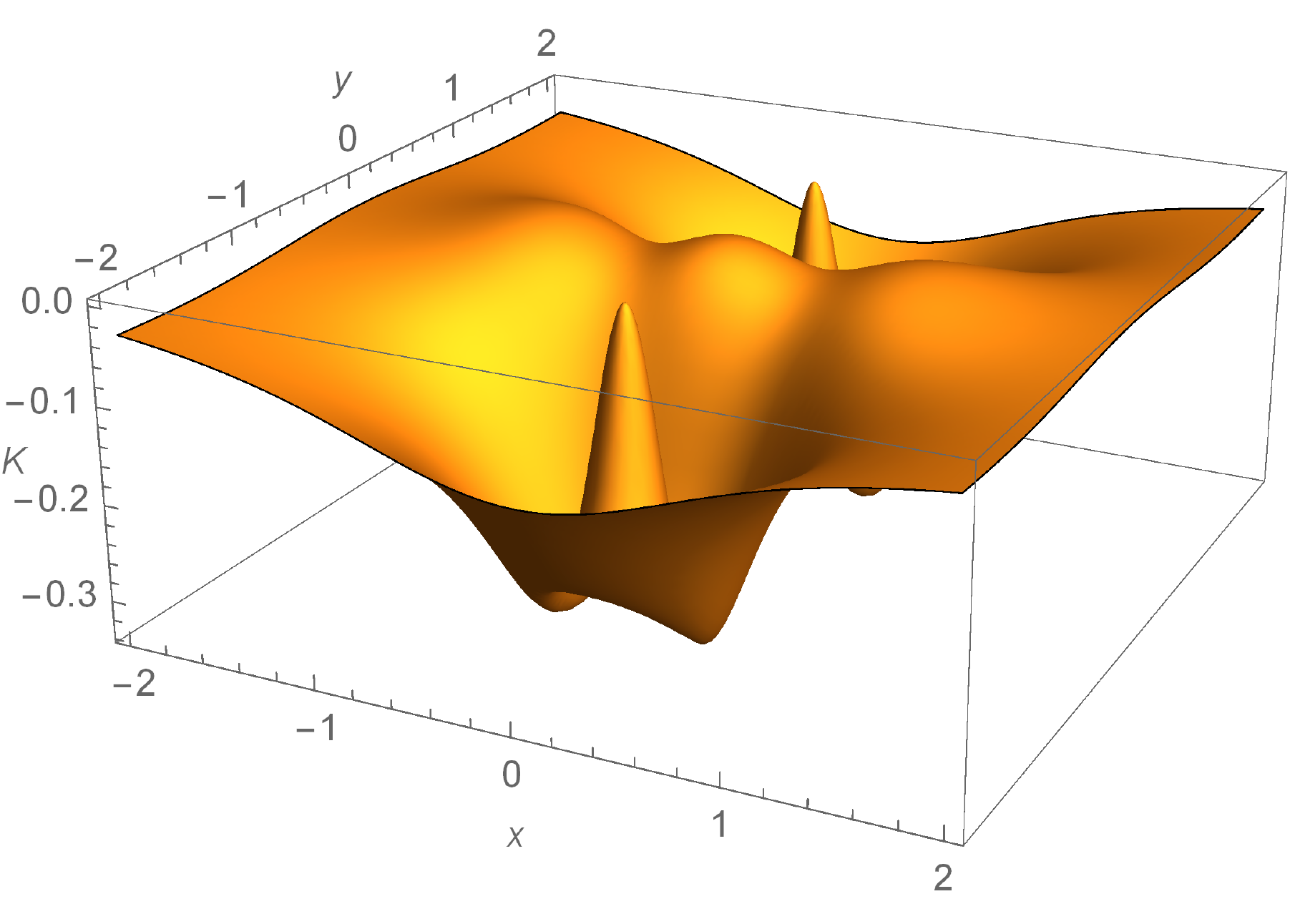}
\end{center}
\caption{Graph of the Gaussian curvature $K(x,y)$ of the surface $\mathcal{P}$ corresponding to configurations that are parallelograms.}
\label{plot:shirt}
\end{figure}

{\sc Remark} As in \cite{Pants}, we consider the Gaussian curvature of $\mathcal{P}$ around the (simultaneous) binary collisions. In the same way as \cite{AGMR}, from equations (\ref{eq:K}) and (\ref{eq:Dlambda}), we obtain
\[
\lim_{(x,y)\rightarrow C}K(x,y)=0,
\]
where $C$ is some (simultaneous) binary collision, see Figure \ref{plot:shirt}. It follows that the ends are asymptotic to Euclidean cylinders.
In fact, if $K_{ij}$ is the radius of the Euclidean cylinder when we approach the $ij$ end,  we have that (see \cite[eq. 3.16]{Pants} for explicit computations)
\[
K_{13}=K_{24}=\frac{1}{\sqrt{2}}, \qquad \text{(binary collisions)},
\]
and
\[
K_{12}=K_{14}=1, \quad \text{(simultaneous binary collisions)}.
\]

\textbf{Proof of Thm. \ref{C}:}

Restricting to the totally geodesic collinear configurations amounts to reducing the metric $$U|_{\R^3}ds_{Eucl}^2$$ by dilations which is equivalent to restricting to the sphere $I=1$.

That is in the usual spherical coordinates ($I=\rho^2$) we wish to compute the curvature, $K_C$,  of $$UI(d\phi^2+\cos^2\phi d\theta^2).$$ We set $u=UI=U|_{S^2}$.

We have the standard Lemma \ref{lem:crmetric} (via say isothermal coordinates): $$u K_C=1-\frac12 \Delta_{S^2} \log u$$ where $\Delta_{S^2}$ is the Laplacian on the sphere.

As we have: $$\Delta_{S^2}\log u=(\Delta \log \rho^2 U)|_{S^2}=(\Delta\log U)|_{S^2}+2$$ where $\Delta$ is the standard Euclidean Laplacian, we find that $K_C< 0$ is equivalent to $(\Delta \log U)|_{S^2}> 0$.

So it suffices to show that $\Delta \log U> 0$ as we do now:

Take the coordinates:

$$2\xi_1=q_1-q_2-q_3+q_4,~~ 2\xi_2=q_1-q_2+q_3-q_4,$$ $$ 2\xi_3=q_1+q_2-q_3-q_4$$

where we have $$U=\sum_{1\le i<j\le 3} \frac{1}{(\xi_i-\xi_j)^2}+\frac{1}{(\xi_i+\xi_j)^2}.$$

And introducing the notations $$s_{ij}\inv=\xi_i+\xi_j,~~~r_{ij}\inv=\xi_i-\xi_j$$ $$^kx_{ij}=s_{ij}^k+r_{ij}^k=(\xi_1+\xi_j)^{-k}+(\xi_i-\xi_j)^{-k},$$ we find (with $\del_i=\frac{\del}{\del\xi_i}$ and for $1\le i,j\le 3$): $$U=\sum_{i<j} ~^2x_{ij}$$ $$\del_i^2U=6\sum_{j\ne i} ~^4x_{ij}$$ $$\del_iU=-2\sum_{j\ne i}~^3x_{ij}.$$

Now by setting $f(x,y)=3x^2+3y^2-4xy=x^2+y^2+2(x-y)^2> 0$ (for $(x,y)\ne (0,0)$) and using $$U^2\del_i^2\log U=U\del_i^2U-(\del_iU)^2$$ one computes: $$U^2\del_1^2\log U=2(^6x_{12}+^6x_{13}+s_{12}^2r_{12}^2f(s_{12},r_{12})+s_{13}^2r_{13}^2f(s_{13},r_{13}))+$$ $$2(s_{12}^2s_{13}^2f(s_{12},s_{13})+s_{12}^2r_{13}^2f(s_{12},r_{13})+s_{13}^2r_{12}^2f(s_{13},r_{12})+r_{12}^2r_{13}^2f(r_{12},r_{13}))+$$ $$^2x_{23}~^4x_{12}+^2x_{23}~^4x_{13}$$ which is a sum of positive terms and likewise for $\del_2^2\log U, \del_3^2\log U$ so that $\Delta\log U > 0$.\qed

\textbf{Proof of Thm \ref{M}:}

After the spherical projection $$u=\frac{\xi_1+\xi_2}{1-\xi_3}, v=\frac{\xi_1-\xi_2}{1-\xi_3}$$ our metric $U(\sum d\xi_i^2)|_{S^2}$ becomes $\frac{\lambda}{2}(du^2+dv^2)$ with:

$$\lambda=\frac{1}{u^2}+\frac{1}{v^2}+\frac{16}{((u-1)^2+(v-1)^2-4)^2}+$$ $$\frac{16}{((u+1)^2+(v+1)^2-4)^2}+\frac{16}{((u+1)^2+(v-1)^2-4)^2}+\frac{16}{((u-1)^2+(v+1)^2-4)^2}.$$\\

\begin{figure}[ht!]
\centering
\begin{tikzpicture}[scale=1.50]
\draw [<->, thick, blue] (0,4) -- (8,4); 
\draw [<->, thick, blue] (4,0) -- (4,8); 
\draw [ thick, blue] (5,5) circle [radius=2];;
\draw [ thick, blue] (3,5) circle [radius=2];;
\draw [ thick, blue] (5,3) circle [radius=2];;
\draw [ thick, blue] (3,3) circle [radius=2];;

\draw [<->,very thick, green] (0,8) -- (8,0);
\draw [<-,very thick, green] (0,0) -- (3.7,3.7);
\draw [->,very thick, green] (3.9,3.9) -- (8,8);
\draw [very thick, green] (4,4) circle [radius=1.414];;


\node [below right] at (3.6,3.9) {$T$};
\node [below] at (8,4) {$u$};

\node [right] at (4,8) {$v$};


\end{tikzpicture}
\caption{~~The region $T$ we are focusing on in $u,v$ coordinates. Binary collisions are blue. The equators of the shirts are shown in green.}
\label{region}
\end{figure}
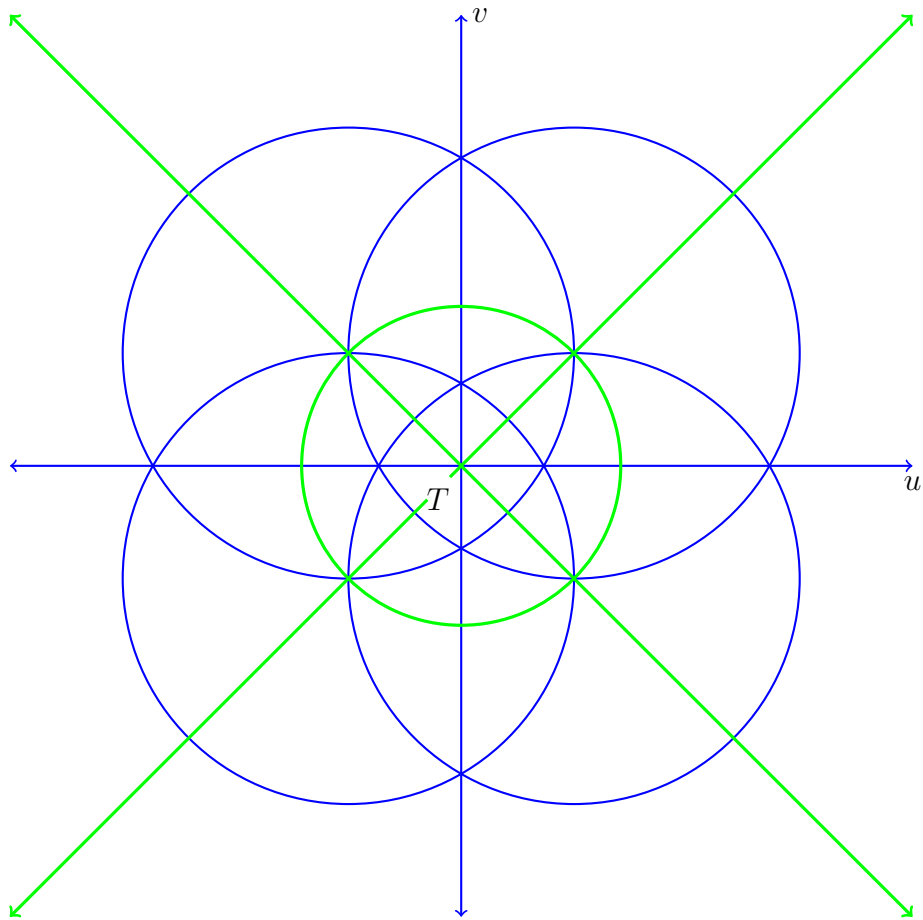

We focus on the invariant region $T$ (see Figure \ref{region}) corresponding to the ordering $q_1\le q_2\le q_3\le q_4$. Here the boundary in $u, v$ coordinates translates to: $$u=0\iff q_1=q_2$$ $$v=0\iff q_3=q_4$$ $$w=(u-1)^2+(v-1)^2-4=0\iff q_2=q_3.$$

Near a binary collision point (for example near $q_1=q_2$ i.e. in a region $0< u<<1$ and $v,w\ge\delta>0$) the metric takes the form
\begin{equation}
    (\frac{1}{2u^2}+O(1))(du^2+dv^2).
    \label{bc}
\end{equation}

Near the simultaneous binary collision, i.e. in the region $0< u,v<<1$, through the change of variable $u+iv=z\mapsto z^2$, the metric takes the same form as eq. \ref{bc}.

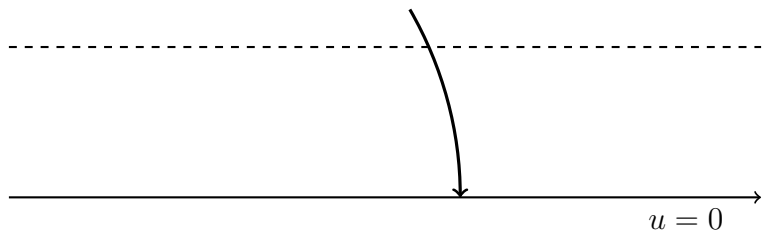
\begin{figure}[ht!]
\centering
\begin{tikzpicture}
\draw [->,thick] (0,2) -- (10,2); 
\draw [dashed,,thick] (0,4) -- (10,4); 
\draw[<-,very thick] (6,2) arc (0:30:5cm);
\node [below] at (9,2) {$u=0$};
\end{tikzpicture}
\caption{~~Sufficiently near the boundary geodesics behave as in the upper half plane. See the Appendix for details.}
\label{hit}
\end{figure}

We fix our attention to unit speed geodesics and note that by eq. \ref{bc} above if a geodesic $\gamma$ passes sufficiently close to the boundary it behaves as a geodesic in the hyperbolic plane, for instance we must have $\alpha(\gamma)\cap\del T$ be a point, see Figure \ref{hit} and the appendix.

Now fix some $q\in \del T\backslash\{$ triple collisions$\}$ and $p\in\del T\backslash\{ q\}$.  We first show the existence of a geodesic from $p$ to $q$.

\begin{figure}[ht!]
\centering
\begin{tikzpicture}
\draw [->,thick] (0,2) -- (11,2); 

\node [above] at (3,5) {$I$};
\draw [dashed,,thick] (6,5) -- (11,5); 
\draw [ultra thick] (0,5) -- (6,5); 
\draw[very thick] (3,2) arc (0:90:3cm);
\draw[very thick] (9,2) arc (0:180:3cm);

\shade[ball color=black] (3,2) circle (.06);
\node [below] at (3,2) {$q$};

\shade[ball color=black] (10,2) circle (.06);
\node [below] at (10,2) {$p$};
\end{tikzpicture}
\caption{~~The compact set $I$. See the Appendix for details.}
\label{comp}
\end{figure}
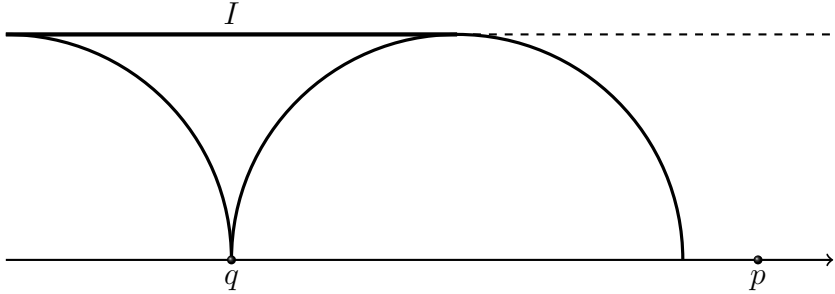

\begin{figure}[ht!]
\centering
\includegraphics[width=50mm] {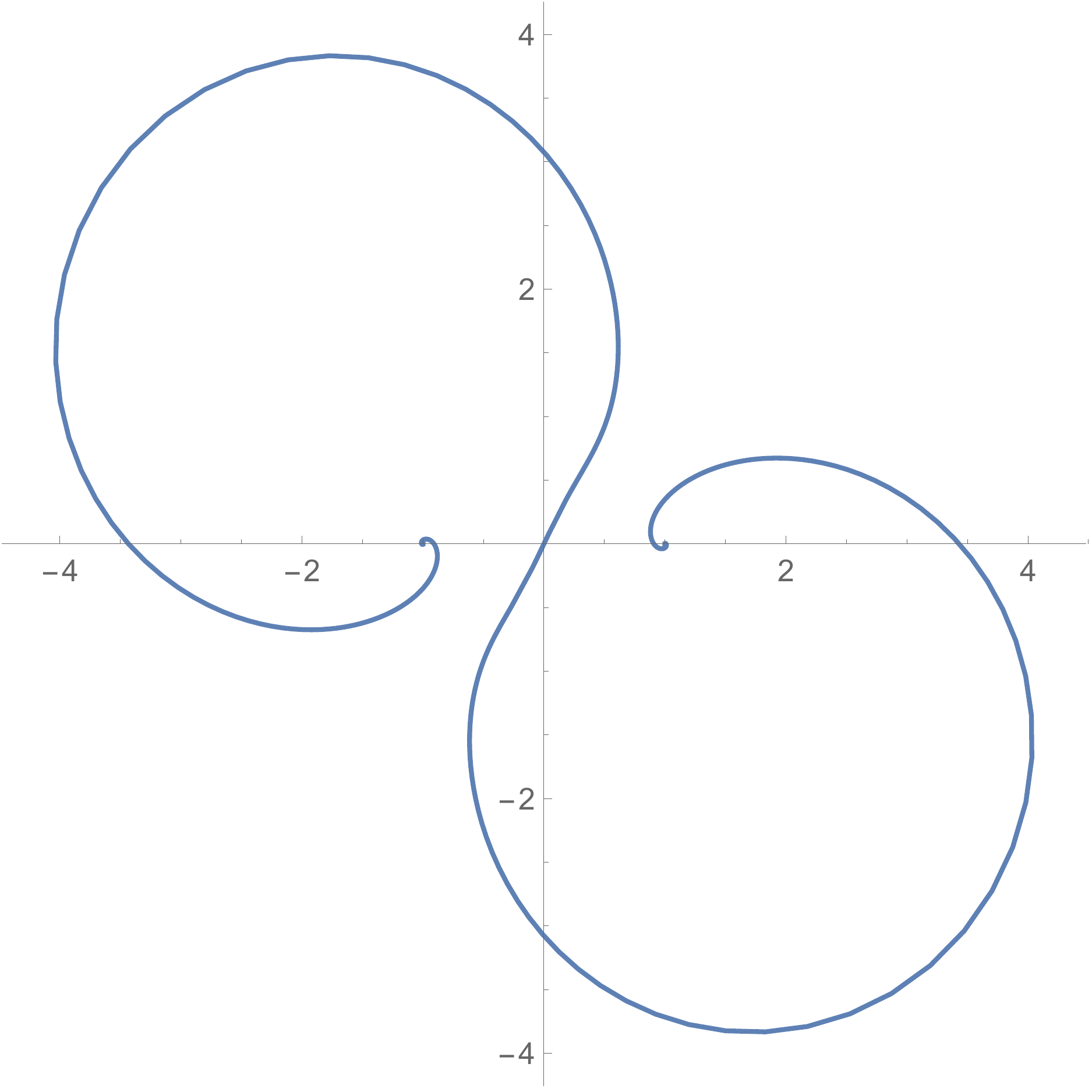}
\includegraphics[width=50mm] {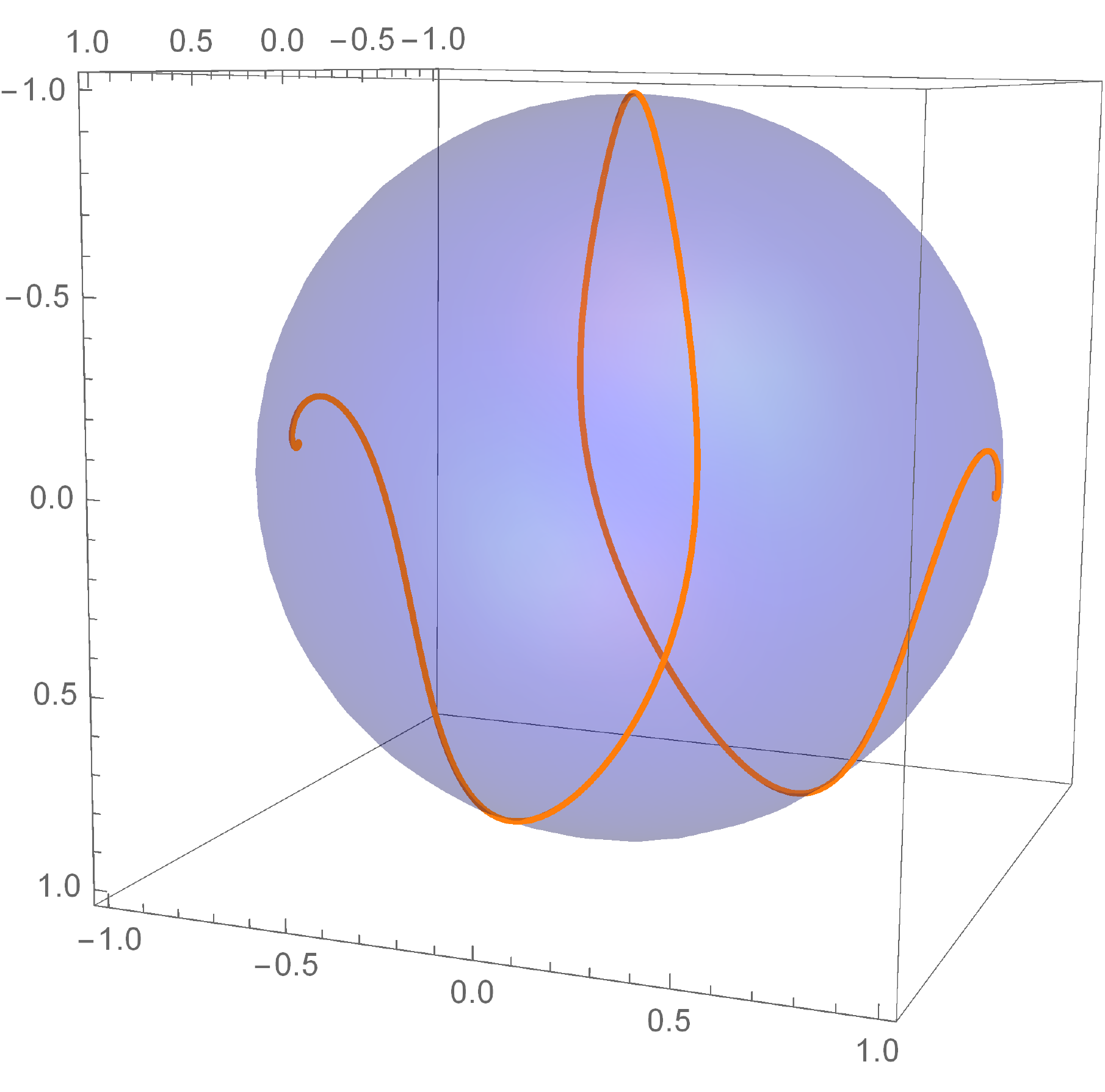}
\caption{~~A distinct collision solution seen from the $(x,y)$-plane ({\em on the left}) and in the shape sphere ({\em on the right}). In this case the collisions are simultaneous binary collisions. }
\label{sol1}
\end{figure}

\begin{figure}[ht!]
\centering
\includegraphics[width=55mm,height =50mm] {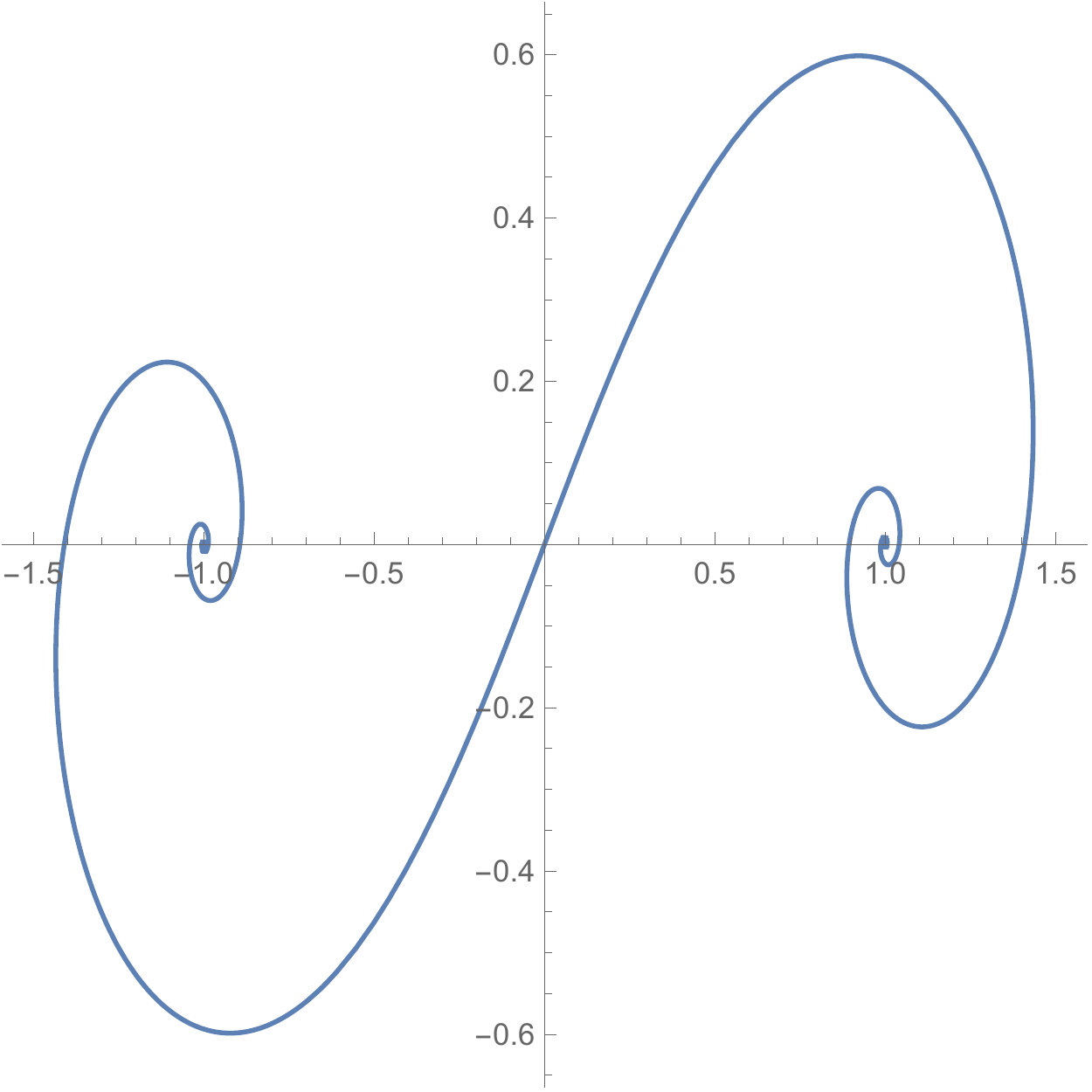}
\includegraphics[width=50mm] {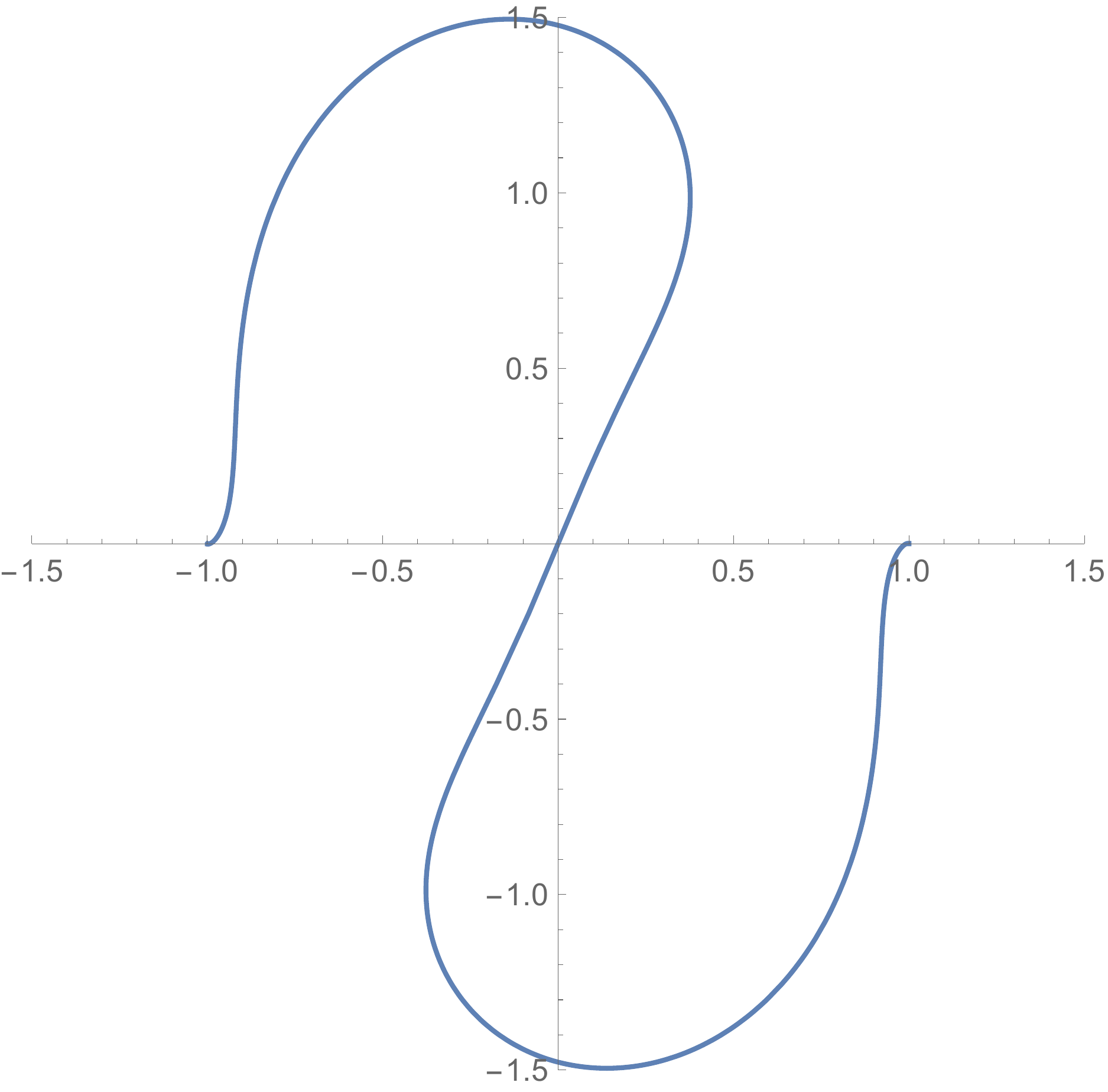}
\includegraphics[width=50mm] {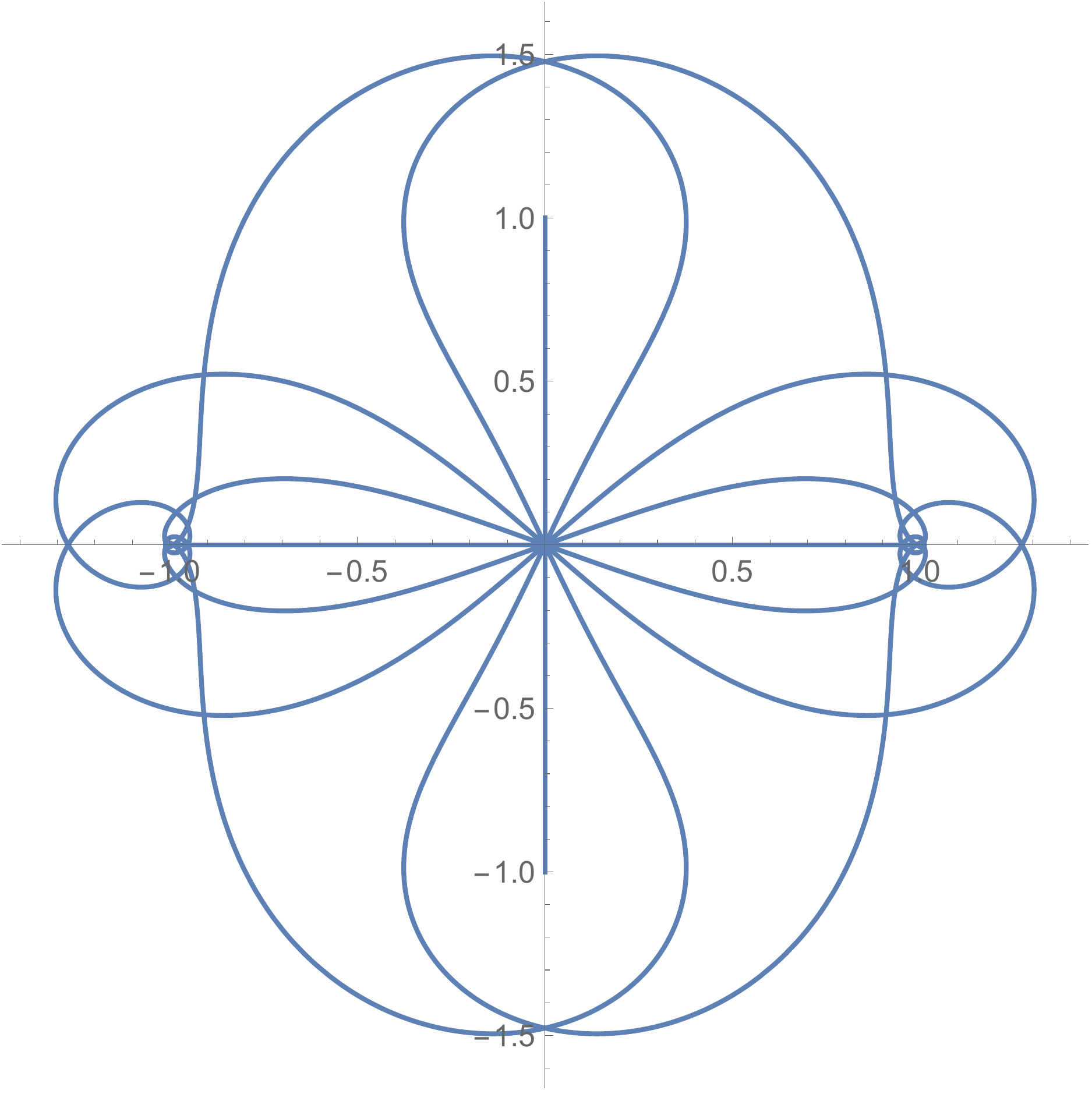}
\caption{~~Solutions with collision seen from the $(x,y)$-plane.
Note that all these solutions pass through a square configuration, corresponding to the south pole in the shape sphere.}
\label{sol2}
\end{figure}

Near $q$  (by eq. \ref{bc}) we have that every geodesic beginning in a neighborhood of $p$ with $\alpha(\gamma)=q$ intersects some compact set $I$, see Figure \ref{comp}.

Now take a sequence of points $q_i\in T^o$ with $q_i\to q$ and a sequence of points $p_i\in T^o$ with $p_i\to p$. Let $\gamma_i$ be the geodesic from $p_i$ to $q_i$. Then for $i$ large enough all such geodesics intersect $I$ in some point $x_i=\gamma_i(t_i)$. Let $v_i=\dot\gamma_i(t_i)$. Now due to the compactness, $(x_i, v_i)\to (x,v)\in I\times \R^2$ and the geodesic with this initial condition $\gamma_{(x,v)}$ goes from $p$ to $q$.

For uniqueness, we note that when we have both $p, q\in \del T\backslash \{$ triple collisions $\}$ then the angles at the boundary between any two geodesics from $p$ to $q$ are zero (see appendix). So if there were two such geodesics bounding a region $R$, the Gauss-Bonnet theorem yields $$2\pi=\int_R K~dA+\int_{\del R}\kappa_g~ds=\int_R KdA+\pi +\pi.$$ Which is impossible by the negative curvature. Hence in this case such a geodesic is unique. \qed\\

{\sc Remark} Due to the covering of the exponential map, every geodesic passes sufficiently close to the boundary and hence is characterized uniquely by the points $\alpha(\gamma), \omega(\gamma)$ (except possibly for geodesics beginning or ending in triple collision). The dynamics here on $\mathcal{C}$ are an attractive case of the \textit{integrable} Calogero-Moser system.

\section{Acknowledgments}
We would like to thank Richard Montgomery for encouragement and discussion regarding perturbed metrics along with Jie Qing. Also we thank Rick Moeckel who posed the question of theorem \ref{M}, Gabriel Martins for his interest and questions on the parallelogram space.

\begin{appendix}
\section{On a metric of the form of eq. \ref{bc}}

Consider geodesics of the metric $$(\frac{1}{y^2}+f(x,y))(dx^2+dy^2)=\frac{\phi}{y^2}(dx^2+dy^2)$$ on the upper half plane $z=x+iy$, $y>0$. Where $\phi,\frac{1}{\phi}\in O(1)$, $\phi_x\in O(y^2)$ and $\phi_y\in O(y)$ as $y\to 0$. We recast in Hamiltonian form, that is with $$H=\frac{y^2}{2\phi}(p_x^2+p_y^2)$$ and the equations of motion for $H=1$ (so $yp_x, yp_y\in O(1)$) read:

$$\dot x=\frac{y^2}{\phi}p_x,  \dot y=\frac{y^2}{\phi}p_y,$$
$$\dot p_x=\frac{\phi_x}{\phi},  \dot p_y=-\frac{2}{y}+\frac{\phi_y}{\phi}.$$

In particular for $y$ sufficiently small we have $\dot p_y<0$ so that if we enter such a sufficiently near region of the boundary ($p_y(0)<0\so p_y(t)<0$) then $\dot y<0$ and $y\to 0$ as $t\to \infty$ (see Figure \ref{hit}). Moreover near the boundary we may parametrize geodesics by $y$ as they approach or leave the boundary.

Parametrizing now by $y$ and letting $'$ denote $\frac{d}{dy}$ we have: $$p_yp_y'=-\frac{2}{y^3}+O(\frac{1}{y})\so p_y^2=\frac{1}{y^2}+O(\log y)+c_y$$ and $$p_x'=O(\frac{1}{p_y})=O(y)\so p_x=O(y^2)+c_x$$ in particular, $p_y\to\infty$ and $p_x\to c_x$ as $y\to 0$ so that $\theta=\arctan(\dot y/\dot x)=\arctan(p_y/p_x)\to \pi/2$ i.e. all orbits intersect the boundary perpendicularly.\\

Take an $\e$ small enough so that the above observations hold over $y\le\e$.\\

Now we direct our attention towards establishing Figure \ref{comp}.

We reparametrize by $y$, to obtain $$|\frac{dx}{dy}|=|\frac{p_x}{p_y}|=\frac{y|p_x|}{\sqrt{1+y^2O(\log(y))+c_y y^2}}=\frac{y|p_x|}{1+O(y)}.$$

In particular as $p_x=O(y^2)+c_x$, any orbit entering $y\le\e$ with an initial $p_x(\e)=\delta$ at the instant $y=\e$ has $$|p_x|\le \delta+O(\e^2)$$ over the rest of it's fall so that $$|x'|\le c_1y(\delta+O(\e^2))$$ and $$|x(\e)-x(y)|\le c_1\e^2\delta+O(\e^4)$$ for some constant $c_1$.

Moreover by energy conservation, $H=1$, we have $\delta=|p_x(\e)|\le \frac{\sqrt{2\phi}}{\e}$, so that $$\delta\le\frac{c_2}{\e}$$ for some constant $c_2$.

Hence any geodesic entering $y\le\e$ at $x(\e)+i\e$ will hit the boundary at a point $x(0)$ with $$|x(\e)-x(0)|\le c\e$$ for some constant $c$. Likewise if we reverse the time a geodesic leaving the region $y\le\e$ at $x(\e)+i\e$ will have come from a boundary point $x(0)=x(\e)+O(\e)$.

Divide those geodesics coming into $q$ into two classes:
\begin{enumerate}
\item[1.] Those that do not intersect $y=\e$,
\item[2.]  Those that intersect $y=\e$.
\end{enumerate}

For the geodesics of class 1, we have by $\dot p_y<0$ that there is a unique maximum attained at $y=M<\e$. Now our estimates above imply we intersect the boundary again a distance at most $2cM<2c\e$ from $q$.

Consider geodesics of class 2 intersecting $y=\e$ at $q+x(\e)+i\e$. Then by the above, once $|x(\e)|>c\e$ such a geodesic cannot reach $q$, hence all such geodesics of class 2 intersect the compact set $\{q+x+iy: -c\e\le x\le c\e, y=\e\}$.

In particular we can establish Figure 8 by choosing $\e$ small enough that $p$ does not lie within $c\e$ of $q$.

\end{appendix}


\begin{thebibliography}{15}

\bibitem{AGMR} M. Alvarez-Ram\'irez, A. Garc\'ia, J. Mel\'endez and J. Guadalupe Reyes-Victoria, \textit{The three-body problem with equal masses via the hyperbolic pants and equivariant Riemannian geometry}, Preprint 2016. 

\bibitem{ACM} V.I. Arnold, \textit{Mathematical Methods of Classical Mechanics}, 2nd ed. 1997, Springer.

\bibitem{Chen} KC. Chen, \textit{Action Minimizing orbits in the Parallelogram Four-Body Problem with Equal Masses}, Archive for Rational Mechanics and Analysis, July 2001, Volume 158, Issue 4, pp 293--318.

\bibitem{Nopants} C. Jackman, R. Montgomery, \textit{No Hyperbolic pants for the 4-body problem with strong potential}, Pacific Journal of Mathematics 280-2 (2016), 401--410.

\bibitem{Pants} R. Montgomery, \textit{Hyperbolic Pants fit a three-body problem}, Erg. Th. and Dyn. Systems, v. 25, (2005), 921--947.

\bibitem{Shape} R. Montgomery, \textit{The three body problem and the shape sphere}, Amer. Math. Monthly, v 122, no. 4, (2015) pp 299--321.


\end{thebibliography}
\end{document}